\documentclass[10pt]{article}
\usepackage{cite}
\usepackage{mathrsfs}
\usepackage{amsfonts}
\usepackage{amsmath}
\usepackage{amsfonts,amssymb}
\usepackage{dsfont}
\usepackage{curves}
\usepackage{mathrsfs}
\usepackage{pifont}
\usepackage{amssymb}
\allowdisplaybreaks

\numberwithin{equation}{section}

\date{}

\def\BigRoman{\uppercase\expandafter{\romannumeral\number\count 255 }}
\def\Romannumeral{\afterassignment\BigRoman\count255=}

\begin{document}
\title{Sufficient conditions for fractional $[a,b]$-deleted graphs
\thanks{This work was supported by the Natural Science Foundation of Shandong Province (ZR2023MA078)}
}
\author{\small  Sizhong Zhou$^{1}$, Yuli Zhang$^{2}$\footnote{Corresponding
author. E-mail address: zhangyuli\_djtu@126.com (Y. Zhang)}\\
\small $1$. School of Science, Jiangsu University of Science and Technology,\\
\small Zhenjiang, Jiangsu 212100, China\\
\small $2$. School of Science, Dalian Jiaotong University,\\
\small Dalian, Liaoning 116028, China\\
}

\maketitle
\begin{abstract}
\noindent Let $a$ and $b$ be two positive integers with $a\leq b$, and let $G$ be a graph with vertex set $V(G)$ and edge set $E(G)$. Let
$h:E(G)\rightarrow[0,1]$ be a function. If $a\leq\sum\limits_{e\in E_G(v)}{h(e)}\leq b$ holds for every $v\in V(G)$, then the subgraph of
$G$ with vertex set $V(G)$ and edge set $F_h$, denoted by $G[F_h]$, is called a fractional $[a,b]$-factor of $G$ with indicator function
$h$, where $E_G(v)$ denotes the set of edges incident with $v$ in $G$ and $F_h=\{e\in E(G):h(e)>0\}$. A graph $G$ is defined as a fractional
$[a,b]$-deleted graph if for any $e\in E(G)$, $G-e$ contains a fractional $[a,b]$-factor. The size, spectral radius and signless Laplacian
spectral radius of $G$ are denoted by $e(G)$, $\rho(G)$ and $q(G)$, respectively. In this paper, we establish a lower bound on the size,
spectral radius and signless Laplacian spectral radius of a graph $G$ to guarantee that $G$ is a fractional $[a,b]$-deleted graph.
\\
\begin{flushleft}

(2020) Mathematics Subject Classification: 05C50, 05C70, 05C72
\end{flushleft}
\end{abstract}

\section{Introduction}

In this paper, we only deal with finite undirected graphs which have neither loops nor multiple edges. Let $G$ be a graph with vertex set
$V(G)=\{v_1,v_2,\cdots,v_n\}$ and edge set $E(G)$. The order and size of $G$ are denoted by $|V(G)|=n$ and $|E(G)|=e(G)$, respectively. For
a vertex $v$ of $G$, the degree of $v$ in $G$, denoted by $d_G(v)$, is the number of vertices in $G$ which are adjacent to $v$. Let
$\delta(G)=\min\{d_G(v):v\in V(G)\}$ denote the minimum degree of $G$. Let $N_G(v)$ denote the neighborhood of a vertex $v$ in $G$. For
any $S\subseteq V(G)$, we use $G[S]$ to denote the subgraph of $G$ induced by $S$, and write $G-S=G[V(G)\setminus S]$. A vertex subset
$S\subseteq V(G)$ is called independent if $G[S]$ has no edges. For two vertices subsets $S,T\subseteq V(G)$ with $S\cap T=\emptyset$, let
$E_G(S,T)$ denote the set of edges admitting one end-vertex in $S$ and the other in $T$ and set $e_G(S,T)=|E_G(S,T)|$. As usual, the complete
graph of order $n$ is denoted by $K_n$. Let $G_1$ and $G_2$ be two vertex disjoint graphs. We use $G_1\cup G_2$ to denote the disjoint union
of $G_1$ and $G_2$. The join $G_1\vee G_2$ is the graph obtained from $G_1\cup G_2$ by adding all the edges joining a vertex of $G_1$ to a
vertex of $G_2$.

Recall that $V(G)=\{v_1,v_2,\cdots,v_n\}$. The adjacency matrix $A(G)=(a_{ij})$ of $G$ is the $n\times n$ symmetric matrix, where $a_{ij}=1$
if $v_i$ and $v_j$ are adjacent in $G$, zero otherwise. Let $D(G)$ be the diagonal degree matrix of $G$. Then $Q(G)=D(G)+A(G)$ is called the
signless Laplacian matrix of $G$. The largest eigenvalues of $A(G)$ and $Q(G)$, denoted by $\rho(G)$ and $q(G)$, are called the spectral radius
and the signless Laplacian spectral radius of $G$, respectively.

Let $g$ and $f$ be two integer-valued functions defined on $V(G)$ with $0\leq g(v)\leq f(v)$ for any $v\in V(G)$. A $(g,f)$-factor of $G$ is
a spanning subgraph $F$ of $G$ such that $g(v)\leq d_F(v)\leq f(v)$ for all $v\in V(G)$. Let $a$ and $b$ are two integers with $1\leq a\leq b$.
If $g\equiv a$ and $f\equiv b$, then a $(g,f)$-factor of $G$ is called an $[a,b]$-factor of $G$. A $[k,k]$-factor of $G$ is called a $k$-factor
of $G$. In particular, a 1-factor is also called a perfect matching. Let $h:E(G)\rightarrow[0,1]$ be a function. If
$g(v)\leq\sum\limits_{e\in E_G(v)}{h(e)}\leq f(v)$ holds for every $v\in V(G)$, then the subgraph of $G$ with vertex set $V(G)$ and edge set
$F_h$, denoted by $G[F_h]$, is called a fractional $(g,f)$-factor of $G$ with indicator function $h$, where $E_G(v)$ denotes the set of edges
incident with $v$ in $G$ and $F_h=\{e\in E(G):h(e)>0\}$. If $g\equiv a$ and $f\equiv b$, then a fractional $(g,f)$-factor of $G$ is called a
fractional $[a,b]$-factor of $G$. In particular, a fractional $(f,f)$-factor of $G$ is called a fractional $f$-factor of $G$.

In mathematical literature, the study on factors and fractional factors of graphs attracted much attention. Some sufficient conditions for
graphs to possess $[1,2]$-factors were derived by Kelmans \cite{Ke}, Kano, Katona and Kir\'{a}ly \cite{KKK}, Zhou, Wu and Bian \cite{ZWB},
Zhou and Bian \cite{ZB}, Zhou, Sun and Liu \cite{ZSL}, Zhou, Wu and Xu \cite{ZWX}, Zhou \cite{Zhr,Zd}, Gao and Wang \cite{GW}, Liu \cite{Ls},
Wang and Zhang \cite{WZi}, Wu \cite{Wp}. Much effort has been devoted to obtain some results on $[a,b]$-factors in graphs by utilizing
various graphic parameters such as Fan-type condition \cite{M}, toughness \cite{Ka}, stability number \cite{KL}, degree condition \cite{WZo}
and others. Gao, Wang and Chen \cite{GWC}, Wang and Zhang \cite{WZr}, Zhou \cite{Za1,Za2,Za3,Zr}, Zhou, Liu and Xu \cite{ZLX} presented some
sufficient conditions for graphs to possess fractional $[a,b]$-factors.

For any function $\varphi$ defined on $V(G)$, write $\varphi(S)=\sum\limits_{v\in S}\varphi(v)$, where $S\subseteq V(G)$. In particular,
$\varphi(\emptyset)=0$. In 1970, Lov\'asz \cite{L} characterized a graph with a $(g,f)$-factor.

\medskip

\noindent{\textbf{Theorem 1.1}} (\cite{L}). Let $G$ be a graph, and let $g,f$ be two nonnegative integer-valued functions defined on $V(G)$
satisfying $g(v)\leq f(v)$ for each $v\in V(G)$. Then $G$ has a $(g,f)$-factor if and only if for all disjoint subsets $S$ and $T$ of $V(G)$,
$$
f(S)+d_{G-S}(T)-g(T)-q_G(S,T,g,f)\geq0,
$$
where $q_G(S,T,g,f)$ denotes the number of components $C$ of $G-(S\cup T)$ satisfying $g(v)=f(v)$ for all $v\in V(C)$ and
$f(V(C))+e_G(C,T)\equiv1 \ (mod \ 2)$.

\medskip

In 1990, Anstee \cite{A} posed a criterion for a graph admitting a fractional $(g,f)$-factor. Liu and Zhang \cite{LZ} presented a new proof.

\medskip

\noindent{\textbf{Theorem 1.2}} (\cite{A,LZ}). Let $G$ be a graph, and let $g,f$ be two functions from $V(G)$ to the nonnegative integers with
$g(v)\leq f(v)$ for every $v\in V(G)$. Then $G$ has a fractional $(g,f)$-factor if and only if
$$
f(S)+d_{G-S}(T)-g(T)\geq0
$$
for any $S\subseteq V(G)$, where $T=\{v:v\in V(G)\setminus S, d_{G-S}(v)<g(v)\}$.

\medskip

A graph $G$ is defined as a fractional $(g,f)$-deleted graph if for any $e\in E(G)$, $G-e$ contains a fractional $(g,f)$-factor. In 2003, Li,
Yan and Zhang \cite{LYZ} provided a criterion for a graph to be a fractional $(g,f)$-deleted graph, which is an extension of Theorem 1.2. For
two integers $a$ and $b$ with $1\leq a\leq b$, if $g\equiv a$ and $f\equiv b$, then a fractional $(g,f)$-deleted graph is called a fractional
$[a,b]$-deleted graph. A fractional $[a,b]$-deleted graph is called a fractional $k$-deleted graph when $a=b=k$. Kotani \cite{Ko} obtained a
binding number condition for the existence of fractional $k$-deleted graphs. The following theorem on fractional $[a,b]$-deleted graph is a
special case of Li, Yan and Zhang's fractional $(g,f)$-covered graph theorem.

\medskip

\noindent{\textbf{Theorem 1.3}} (\cite{LYZ}). Let $G$ be a graph and $a\leq b$ be two positive integers. Then $G$ is a fractional
$[a,b]$-deleted graph if and only if
$$
b|S|+d_{G-S}(T)-a|T|\geq\varepsilon(S,T)
$$
for any $S\subseteq V(G)$, where $T=\{v:v\in V(G)\setminus S, d_{G-S}(v)\leq a\}$ and $\varepsilon(S,T)$ is defined by
\[
 \varepsilon(S,T)=\left\{
\begin{array}{ll}
2,&if \ T \ is \ not \ independent,\\
1,&if \ T \ is \ independent \ and \ e_G(T,V(G)\setminus(S\cup T))\geq1,\\
0,&otherwise.\\
\end{array}
\right.
\]

\medskip

Very recently, O \cite{Os} derived a spectral radius condition to guarantee that a graph contains a 1-factor (or a perfect matching). Zhou and
Liu \cite{ZL} presented a spectral radius condition for the existence of an odd $[1,b]$-factor in a graph. Motivated by \cite{Os,ZL,LYZ} directly,
it is natural and interesting to show some spectral radius conditions to guarantee the existence of fractional $[a,b]$-deleted graphs. In what
follows, we put forward a spectral radius condition and a signless Laplacian spectral radius condition for a graph to be a fractional
$[a,b]$-deleted graph, respectively.

\medskip

\noindent{\textbf{Theorem 1.4.}} Let $a$ and $b$ be two positive integers with $b\geq\max\{a,3\}$, and let $G$ be a graph of order $n$ with
$n\geq\max\{a+2,7\}$. If
$$
\rho(G)>\rho(K_a\vee(K_{n-a-1}\cup K_1)),
$$
then $G$ is a fractional $[a,b]$-deleted graph.

\medskip

\noindent{\textbf{Remark 1.5.}} In what follows, we exhibit that the condition $\rho(G)>\rho(K_a\vee(K_{n-a-1}\cup K_1))$ declared in Theorem 1.4
is sharp, namely, it cannot be replaced by $\rho(G)\geq\rho(K_a\vee(K_{n-a-1}\cup K_1))$.

Let $a$ and $b$ be two positive integers with $a\leq b$ and $G=K_a\vee(K_{n-a-1}\cup K_1)$. Then $\rho(G)=\rho(K_a\vee(K_{n-a-1}\cup K_1))$.
Set $S=\emptyset$ and $T=V(K_1)$. It is obvious that $T$ is independent and $e_G(T,V(G)\setminus(S\cup T))=a\geq1$. According to the definition
of $\varepsilon(S,T)$, we possess $\varepsilon(S,T)=1$. Thus, we deduce
$$
\theta_G(S,T)=b|S|+d_{G-S}(T)-a|T|=0+a-a=0<1=\varepsilon(S,T).
$$
By virtue of Theorem 1.3, $G$ is not a fractional $[a,b]$-deleted graph.

\medskip

\noindent{\textbf{Theorem 1.6.}} Let $a$ and $b$ be two positive integers with $b\geq\max\{a,3\}$, and let $G$ be a connected graph of order $n$
with $n\geq\max\{a+2,7\}$. If
$$
q(G)>2n-4+\frac{a+1}{n-1}
$$
and
$$
q(G)>q(K_a\vee(K_{n-a-1}\cup K_1)),
$$
then $G$ is a fractional $[a,b]$-deleted graph.

\medskip

\noindent{\textbf{Remark 1.7.}} Next, we claim that the condition $q(G)>q(K_a\vee(K_{n-a-1}\cup K_1))$ declared in Theorem 1.6 cannot be replaced
by $q(G)\geq q(K_a\vee(K_{n-a-1}\cup K_1))$. Let $a$ and $b$ be two positive integers with $a\leq b$ and $G=K_a\vee(K_{n-a-1}\cup K_1)$. Then $q(G)=q(K_a\vee(K_{n-a-1}\cup K_1))$. From Remark 1.5, $G$ is not a fractional $[a,b]$-deleted graph. But, I do not know whether the condition
$q(G)>2n-4+\frac{a+1}{n-1}$ declared in Theorem 1.6 is sharp or not.

\medskip

Finally, we pose a sufficient condition to guarantee a graph $G$ to be a fractional $[a,b]$-deleted graph with respect to its size, which
plays a key role to verify Theorems 1.4 and 1.6.

\medskip

\noindent{\textbf{Theorem 1.8.}} Let $a$ and $b$ be two positive integers with $b\geq\max\{a,3\}$, and let $G$ be a graph of order $n$ with
$n\geq\max\{a+2,7\}$. If $\delta(G)\geq a+1$ and
$$
e(G)\geq\binom{n-1}{2}+\frac{a+2}{2},
$$
then $G$ is a fractional $[a,b]$-deleted graph.

\medskip

The proofs of Theorems 1.4 and 1.6 will be provided in Sections 3 and 4, respectively. The proof of Theorem 1.8 will be provided in Section 2.

\section{The proof of Theorem 1.8}

In this section, we prove Theorem 1.8, which gives a sufficient condition to ensure that a graph $G$ is a fractional $[a,b]$-deleted graph
by utilizing Theorem 1.3.

\medskip

\medskip

\noindent{\it Proof of Theorem 1.8.} Let $\theta_G(S,T)=b|S|+d_{G-S}(T)-a|T|$ for any $S\subseteq V(G)$ and
$T=\{v:v\in V(G)\setminus S, d_{G-S}(v)\leq a\}$. Suppose, to the contrary, that $G$ is not a fractional $[a,b]$-deleted graph. According to
Theorem 1.3, there exists some subset $S$ of $V(G)$ such that
\begin{align}\label{eq:2.1}
\theta_G(S,T)=b|S|+d_{G-S}(T)-a|T|\leq\varepsilon(S,T)-1,
\end{align}
where $T=\{v:v\in V(G)\setminus S, d_{G-S}(v)\leq a\}$.

\medskip

\noindent{\bf Claim 1.} $n\geq a+3$.

\noindent{\it Proof.} Let $n=a+2$. In terms of $\delta(G)\geq a+1$, we infer that $G$ is a complete graph of order $n=a+2$. Consequently, for
any $e\in E(G)$, $G$ has a Hamiltonian cycle $C$ such that $e\in E(C)$. Thus, $G-E(C)$ is a $a$-factor of $G$, and also a fractional $a$-factor
of $G$. Obviously, $G$ is a fractional $[a,b]$-deleted graph, which is contrary to the assumption. Hence, $n\geq a+3$. This completes the proof
of Claim 1. \hfill $\Box$

\medskip

\noindent{\bf Claim 2.} $|T|\geq b+1$.

\noindent{\it Proof.} Assume that $|T|\leq b$. Note that $\varepsilon(S,T)\leq|T|$ by the definition of $\varepsilon(S,T)$. Together with
$\delta(G)\geq a+1$, we deduce
\begin{align*}
\theta_G(S,T)=&b|S|+d_{G-S}(T)-a|T|\\
=&b|S|+d_G(T)-e_G(S,T)-a|T|\\
\geq&b|S|+\delta(G)|T|-|S|\cdot|T|-a|T|\\
\geq&b|S|+(a+1)|T|-|S|\cdot|T|-a|T|\\
=&(b-|T|)|S|+|T|\\
\geq&|T|\\
\geq&\varepsilon(S,T),
\end{align*}
which contradicts \eqref{eq:2.1}. Hence, $|T|\geq b+1$. Claim 2 is proved. \hfill $\Box$

In view of Claim 2, we derive
\begin{align}\label{eq:2.2}
n\geq|S|+|T|\geq|S|+b+1.
\end{align}

According to $e(G)\geq\binom{n-1}{2}+\frac{a+2}{2}$, there exist at most $n-1-\frac{a+2}{2}$ edges not in $E_G[V(G)\setminus(S\cup T),T]\cup E(G[T])$.
Thus, we obtain
\begin{align}\label{eq:2.3}
d_{G-S}(T)\geq(n-1-|S|)|T|-2\left(n-1-\frac{a+2}{2}\right).
\end{align}

It follows from \eqref{eq:2.1}, \eqref{eq:2.2}, \eqref{eq:2.3}, $\varepsilon(S,T)\leq2$, $b\geq\max\{a,3\}$, $n\geq\max\{a+2,7\}$, Claims 1 and 2 that
\begin{align*}
\varepsilon(S,T)-1\geq&\theta_G(S,T)=b|S|+d_{G-S}(T)-a|T|\\
\geq&b|S|+(n-1-|S|)|T|-2\left(n-1-\frac{a+2}{2}\right)-a|T|\\
=&b|S|+(n-1-|S|-a)|T|-2\left(n-1-\frac{a+2}{2}\right)\\
\geq&b|S|+(n-1-|S|-a)(b+1)-2\left(n-1-\frac{a+2}{2}\right)\\
=&(b-3)n+(n-|S|-b)+n-ab+3\\
\geq&(b-3)(a+3)+1+n-ab+3\\
=&n+3(b-a)-5\\
\geq&n-5\\
\geq&2\\
\geq&\varepsilon(S,T),
\end{align*}
which is a contradiction. This completes the proof of Theorem 1.8. \hfill $\Box$

\section{The proof of Theorem 1.4}

In this section, we first introduce some necessary preliminary results, which will be used to verify Theorem 1.4. Hong, Shu and Fang \cite{HSF},
Nikiforov \cite{N} presented an important upper bound on the spectral radius $\rho(G)$.

\medskip

\noindent{\textbf{Lemma 3.1}} (\cite{HSF,N}). Let $G$ be a graph of order $n$ with minimum degree $\delta(G)$. Then
$$
\rho(G)\leq\frac{\delta(G)-1}{2}+\sqrt{2e(G)-n\delta(G)+\frac{(\delta(G)+1)^{2}}{4}}.
$$

\medskip

The following observation is very useful when we utilize the above upper bound on $\rho(G)$.

\medskip

\noindent{\textbf{Proposition 3.2}} (\cite{HSF,N}). For a graph $G$ of order $n$ with $e(G)\leq\binom{n}{2}$, the function
$$
f(x)=\frac{x-1}{2}+\sqrt{2e(G)-nx+\frac{(x+1)^{2}}{4}}
$$
is decreasing with respect to $x$ for $0\leq x\leq n-1$.

\medskip

Let $A=(a_{ij})_{n\times n}$ and $B=(b_{ij})_{n\times n}$. Define $A\leq B$ if for any $1\leq i,j\leq n$, $a_{ij}\leq b_{ij}$ and $A<B$ if
$A\leq B$ and $A\neq B$.

\medskip

\noindent{\textbf{Lemma 3.3}} (\cite{BP,HJ}). Let $O$ be an $n\times n$ zero matrix, $A=(a_{ij})$ and $B=(b_{ij})$ be two $n\times n$ matrices
with the spectral radius $\rho(A)$ and $\rho(B)$, respectively. If $O\leq A\leq B$, then $\rho(A)\leq\rho(B)$. Furthermore, if $B$ is irreducible
and $O\leq A<B$, then $\rho(A)<\rho(B)$.

\medskip

\medskip

\noindent{\it Proof of Theorem 1.4.} We first verify the following claim.

\medskip

\noindent{\bf Claim 3.} $\delta(G)\geq a+1$.

\noindent{\it Proof.} Assume that $\delta(G)\leq a$. Then there exists a vertex $v\in V(G)$ such that $d_G(v)\leq a$, which implies that
$G\subseteq K_a\vee(K_{n-a-1}\cup K_1)$. According to Lemma 3.3, we infer
$$
\rho(G)\leq\rho(K_a\vee(K_{n-a-1}\cup K_1)),
$$
which contradicts the condition that $\rho(G)>\rho(K_a\vee(K_{n-a-1}\cup K_1))$. Hence, $\delta(G)\geq a+1$. Claim 3 is proved. \hfill $\Box$

\medskip

In terms of Claim 3, Lemma 3.1 and Proposition 3.2, we obtain
\begin{align}\label{eq:3.1}
\rho(G)\leq&\frac{\delta(G)-1}{2}+\sqrt{2e(G)-n\delta(G)+\frac{(\delta(G)+1)^{2}}{4}}\nonumber\\
\leq&\frac{a}{2}+\sqrt{2e(G)-n(a+1)+\frac{(a+2)^{2}}{4}}.
\end{align}

Note that the graph $K_{n-1}$ is a proper subgraph of the graph $K_a\vee(K_{n-a-1}\cup K_1)$ and the adjacency matrices of connected graphs
are irreducible. Then by Lemma 3.3, we have
\begin{align}\label{eq:3.2}
\rho(G)>\rho(K_a\vee(K_{n-a-1}\cup K_1))>\rho(K_{n-1})=n-2.
\end{align}
It follows from \eqref{eq:3.1} and \eqref{eq:3.2} that
$$
e(G)>\frac{(n-1)(n-2)}{2}+\frac{a+1}{2}=\binom{n-1}{2}+\frac{a+1}{2}.
$$
By virtue of the integrity of $e(G)$, we get
\begin{align}\label{eq:3.3}
e(G)\geq\binom{n-1}{2}+\frac{a+2}{2}.
\end{align}

By \eqref{eq:3.3}, Claim 3 and Theorem 1.8, we see that $G$ is a fractional $[a,b]$-deleted graph. This completes the proof of Theorem 1.4. \hfill $\Box$

\section{The proof of Theorem 1.6}

In this section, we first present a necessary preliminary result, which will be used to prove Theorem 1.6. Feng and Yu \cite{FY} verified an upper bound
on $q(G)$, which has been widely used in the literature.

\medskip

\noindent{\textbf{Lemma 4.1}} (\cite{FY}). Let $G$ be a connected graph with $n$ vertices and $e(G)$ edges. Then
$$
q(G)\leq\frac{2e(G)}{n-1}+n-2.
$$

\medskip

\medskip

\noindent{\it Proof of Theorem 1.6.} We first prove the following claim.

\medskip

\noindent{\bf Claim 4.} $\delta(G)\geq a+1$.

\noindent{\it Proof.} Suppose, to the contrary, that $\delta(G)\leq a$. Then there exists a vertex $v\in V(G)$ satisfying $d_G(v)\leq a$, which
yields that $G\subseteq K_a\vee(K_{n-a-1}\cup K_1)$. In terms of Lemma 3.3, we deduce
$$
q(G)\leq q(K_a\vee(K_{n-a-1}\cup K_1)),
$$
which contradicts the condition that $q(G)>q(K_a\vee(K_{n-a-1}\cup K_1)$. Consequently, $\delta(G)\geq a+1$. This completes the proof of Claim 4. \hfill $\Box$

\medskip

In view of Lemma 4.1 and the condition that $q(G)>2n-4+\frac{a+1}{n-1}$, we get
$$
2n-4+\frac{a+1}{n-1}<q(G)\leq\frac{2e(G)}{n-1}+n-2,
$$
which yields that
$$
e(G)>\frac{(n-1)(n-2)}{2}+\frac{a+1}{2}=\binom{n-1}{2}+\frac{a+1}{2}.
$$
According to the integrity of $e(G)$, we obtain
\begin{align}\label{eq:4.1}
e(G)\geq\binom{n-1}{2}+\frac{a+2}{2}.
\end{align}
It follows from \eqref{eq:4.1}, Claim 4 and Theorem 1.8 that $G$ is a fractional $[a,b]$-deleted graph. This completes the proof of Theorem
1.6. \hfill $\Box$

\medskip

\section*{Data availability statement}

My manuscript has no associated data.

\section*{Declaration of competing interest}

The authors declare that they have no conflicts of interest to this work.



\end{document}